\documentclass[12]{article}
\usepackage{amsmath}
\usepackage{amsfonts}
\usepackage{amssymb}
\usepackage{amsthm}
\usepackage{color}
\usepackage{setspace}

\newtheorem{thm}{Theorem}[section]

\newtheorem{de}{Definition}[section]
\newtheorem{exm}{Example}[section]

\newcommand{\go}[1]{\mathfrak{#1}}

\newcommand{\R}{{\rm I}\kern-0.18em{\rm R}}
\newcommand{\1}{{\rm 1}\kern-0.25em{\rm I}}
\newcommand{\E}{{\rm I}\kern-0.18em{\rm E}}
\newcommand{\p}{{\rm I}\kern-0.18em{\rm P}}

\usepackage{epsfig}

\usepackage{fancyhdr}

\usepackage{graphics}

\usepackage{amsopn}  

\usepackage{amssymb}
\usepackage{amsxtra}
\usepackage{amsfonts}

\usepackage{afterpage}
\usepackage{eucal}
\usepackage{eufrak}

\usepackage{enumerate}

\allowbreak

\author{Lev B Klebanov\footnote{e-mail: levkleb@yahoo.com},  Charles University, MFF\\
and\\ Svetlozar T. Rachev, Stony Brook University }
\title{$\nu$-Generalized Hyperbolic Distributions}
\date{}
\begin{document}
\maketitle

\begin{abstract}
A new class of probability distributions closely connected to generalized hyperbolic distributions is introduced. It is more adapted to study the distributions of sums of random number of random variables. The properties of these distributions are studied. It seems, that this class may be useful for asset returns modeling. 
\end{abstract}

\section{Introduction and Definitions}
The family of generalized hyperbolic (GH) distributions is well-known in Probability. It was introduced by Ole Barndorff-Nielsen, who studied it in the context of physics of wind-blown sand \cite{B-N}. These distributions are given by their probability density functions, and their properties are studied rather well. Particularly, GH distributions are infinitely divisible \cite{B-NH}.  Connections of GH with other distributions is given in \cite{Paol}. One of more interesting applications of GH distributions is the field of Finance \cite{EK, GHM}.

On the other hand, there exist a large literature on financial applications of stable and geometric stable  distributions in Finance (see, for example, \cite{GHM}). Geometric stable distributions were introduced in \cite{KMM} and are analogues of classical stable distributions for the case of sums of a random number of random variables, where the number of summands has geometric distribution. The case of general distribution of the number of summands was studied in \cite{KR} (see also \cite{KKR}). Because the number of transactions on financial market is random, it is natural to consider sums of random number of random variables  as an element of a model for asset returns distribution. Our aim here is to introduce a variant of GH distributions, connected to sums of a random number of random variables. We will call corresponding distributions as $\nu$-GH.

Let us give precise definitions.
\begin{de}\label{ng} Let $X,X_1,\ldots ,X_n, \ldots$ be a sequence of independent identically distributed (i.i.d.) random variables. Suppose that $\{\nu_p, \; p\in \Delta\}$ is a family of random variables, independent of the sequence above, taking positive integer values, and such that $\E \nu_p = 1/p$ for all $p \in \Delta$.  We shall say that the random variable $X$ has $\nu$-strictly Gaussian (or $\nu$-strictly normal)
distribution if 
\begin{verse}
	 a) $ X \stackrel{d}{=}p^{1/2}\sum_{j=1}^{\nu_p}X_j$ for all $p \in \Delta$,\\
	 and\\
	 b) the random variable $X$ has finite second moment: $\E X^2< \infty$.
\end{verse}
\end{de}
\begin{de}\label{ns} Suppose that the sequence of random variables $X,X_1,\ldots ,X_n, \ldots $ and the family  $\{\nu_p, \; p\in \Delta\}$ are the same as in Definition \ref{ng}. We shall say that the random variable $X$ has $\nu$-strictly stable distribution with index $\alpha \in (0,2)$ if \\ $ X \stackrel{d}{=}p^{1/\alpha}\sum_{j=1}^{\nu_p}X_j$ for all $p \in \Delta$.
\end{de}
Let ${\mathcal P}_p(z), \; p \in \Delta$ be a family of probability generating functions (p.g.f.) of random variables $\nu_p$. In \cite{KR} there was shown that $\nu$-strictly Gaussian distribution exists if and only if the semigroup  $\go P$, generated by the family ${\mathcal P}_p(z), \; p \in \Delta$ with superposition operation is commutative. In the case of commutativity ${\go P}$ the system of functional equations
\begin{equation}\label{eq1}
\varphi (t) = {\mathcal P}_p(\varphi (p t)), \; \forall p \in \Delta, \; t>0 
\end{equation}
has a solution satisfying initial values
\begin{equation}\label{eq2}
\varphi (0) = 1, \;\; \varphi^{\prime}(0)= -1. 
\end{equation}
This solution is unique, and may be represented as Laplace transform of a distribution function $\mathcal A$ defined on non-negative semiaxes:
\begin{equation}\label{eq3}
\varphi(t) = \int_{0}^{\infty}\exp(-tx)d{\mathcal A}(x).
\end{equation}
We call this function standard solution of Poincare equation. In this case characteristic function of non-degenerate $\nu$-strictly Gaussian distribution has form $f(t)=\varphi(a t^2)$, where $a>0$ (see \cite{KR, KKR}).  

The function $\varphi$ also may be used to define one-to-one map from the set of all infinitely divisible distributions onto the set of all $\nu$-infinite stable distributions. Namely, the following Theorem holds (see \cite{KR, KKR}).
\begin{thm}\label{isom}
If $\varphi$ is standard solution of Poincare equation, then the function
\begin{equation}\label{eq4}
g(t)= \varphi (-\log(f(t))), \; \; t\in \R^1,
\end{equation}
is $\nu$-infinitely divisible characteristic function if and only if characteristic function $f(t)$ is infinite divisible in classical sense.
\end{thm}
Theorem \ref{isom} is a main tool for definition of $\nu$-GH distributions. Let us also note, that all previous definitions as well as Theorem \ref{isom} remains true for random vectors taking values in Euclidean space $\R^d$. Many examples of $\nu$-stable distributions are given in \cite{KKRT}.

\section{$\nu$-GH distrbutions}
\setcounter{equation}{0}
As it was mentioned in previous section, GH distributions are infinitely divisible. Therefore, we may use their characteristic functions instead of the function $f(t)$ in Theorem \ref{isom}. The function $g(t)$ obtained in such way is a natural analog of GH characteristic function for the case of sums of a random number of random summands. 

Let us give some examples.

\begin{exm} \label{ex1} Let $\{\nu_p, \; p \in (0,1)\}$ be a family of random variables with geometric distribution:
\[ \p\{ \nu_p =k\} = p (1-p)^{k-1}, \; k=1,2, \ldots  \]
It is known (see \cite{KMM, KKR}) that in this case standard solution of Poincare equation has the form
\[ \varphi (t) =\frac{1}{1+t}. \]
For one-dimensional case characteristic function of GH distribution has the form
\[ f(t)=\frac{e^{i t \mu } \left(\sqrt{\alpha ^2-\beta ^2} \delta \right)^{\lambda } K_{\lambda}\left[\sqrt{\alpha ^2-(i t+\beta )^2} \delta \right]}{ \left(\sqrt{\alpha ^2-(i t+\beta )^2} \delta \right)^{\lambda } K_{\lambda}\left[\sqrt{\alpha ^2-\beta ^2} \delta \right]}, \]
where $K_{\lambda}$ is modified Bessel function of second kind. Substituting this into equation (\ref{eq4}) we find that the function
\[  g(t)= 1/\Bigl(1-\log[\frac{e^{i t \mu } \left(\sqrt{\alpha ^2-\beta ^2} \delta \right)^{\lambda }  K_{\lambda}\left[\sqrt{\alpha ^2-(i t+\beta )^2} \delta \right]}{\left(\sqrt{\alpha ^2-(i t+\beta )^2} \delta \right)^{\lambda }K_{\lambda}\left[\sqrt{\alpha ^2-\beta ^2} \delta \right]}]\Bigr)\]
is characteristic function of geometric GH distribution (or geo-GH distribution). 
\end{exm}
As it was noted above, it is possible to substitute multivariate GH characteristic function $f(t)$ into (\ref{eq4}) and obtain characteristic function of multivariate geo-GH distribution. We will not do this substitution in explicit form. 

If our GH distribution does not coincide with Student $t$-distribution, then its characteristic function is analytic in a strip, containing real line. This means, this distribution has exponential tails. It is clear, that in this situation, the function $g(t)$ is analytic in a strep, too. Therefore, the tails of geo-GH distribution are exponential, too.  If our GH distribution coincide with Student $t$-distribution, it has only a finite number of moments. It is clear, that the same is true for geo-Student distribution with characteristic function  $g(t)$.

\begin{exm}\label{ex2} Let us consider the following family of random variables $\{\nu_p, \; p \in (1/n^2,\; n=1,2, \ldots)\}$ with probability generating function 
\begin{equation}\label{eq5}
{\mathcal P}_p(z) =\frac{1}{T_{1/\sqrt{p}}}(1/z), 
\end{equation}
where $T_n(x)$ is Chebyshev polynomial of the first kind. In \cite{KKRT} it was proven, that  (\ref{eq5}) really define a family of probability generation functions. Standard solution of Poincare equation has the form (see, \cite{KKRT})
\[ \varphi(t) = \frac{1}{\cosh (\sqrt{2t})}. \]
Using this function we obtain corresponding representation for $g(t)$:
\[ g(t)= \]
\[\sec\left[\sqrt{2}\log^{1/2}\left[\frac{e^{i t \mu } \left(\sqrt{\alpha ^2+(t-i \beta )^2} \delta \right)^{-\lambda } \left(\sqrt{(\alpha
-\beta ) (\alpha +\beta )} \delta \right)^{\lambda } K_{\lambda}\left[\sqrt{\alpha ^2+(t-i \beta )^2} \delta \right]}{K_{\lambda}\left[\sqrt{\alpha ^2-\beta ^2} \delta \right]}\right]\right] \]
\end{exm}
We call distribution with this characteristic function Chebyshev-GH distribution. Tails behavior of Chebyshev-GH distribution is the same as that for GH distribution.


\begin{thebibliography}{500}
 \bibitem{B-N}
 Barndorff-Nielsen, Ole (1977)
 \newblock{Exponentially decreasing distributions for the logarithm of particle size.}
 \newblock{ Proceedings of the Royal Society of London. Series A, Mathematical and Physical Sciences (The Royal Society) 353 (1674): 401–409.}
  
\bibitem{B-NH}
O. Barndorff-Nielsen and Christian Halgreen (1977)
\newblock{Infinite Divisibility of the Hyperbolic and Generalized Inverse Gaussian Distributions,}
\newblock{Zeitschrift für Wahrscheinlichkeitstheorie und verwandte Gebiete.}

\bibitem{Paol}
Marc S. Paolella (2007)
\newblock{Intermediate Probability, A Computational Approach}
\newblock{John Wiley \& Sons Ltd, The Atrium, Southern Gate, Chichester, England}

\bibitem{EK}
Eberlein, E and U. Keller, (1995)
\newblock{Hyperbolic distributions in finance,}
\newblock{ Bernoulli 1, 281-299.}

\bibitem{GHM}
James E. Gentle, Wolfgang Karl Härdle, Yuichi Mori (Editors) (2012)
\newblock{Handbook of Computational Statistics; Concepts and Methods,}
\newblock{Second revised and updated Edition,}
\newblock{Springer-Verlag, Berlin, Heidelberg.}

\bibitem{KMM} 
Klebanov, L. B.,  Manija, G. M., Melamed, I. A. (1987)
\newblock{$\nu_p$-strictly stable laws and estimation of their parameters,}
\newblock{Lecture Notes in Mathematics (Springer-Verlag), Vol.1233, 23-31.}

\bibitem{KR}
Klebanov, L.B., Rachev, S.T. (1996)
\newblock{Sums of a random number of random variables an their approximations with $\nu$- accompanying infinitely divisible laws,}
\newblock{Serdica Math. J., v. 22, No 4, 471-496.}

\bibitem{KKR}
Klebanov L., Kozubowski T.J., Rachev S.T. (2006)
\newblock{Ill-Posed Problems in Probability and Stability of Random Sums, }
\newblock{Nova Scientific Publishers, New York.}

\bibitem{KKRT}
Klebanov L.B., Kakosyan A.V., Rachev S.T., Temnov G. (2012)
\newblock{On a class of distributions stable under random summation,}
\newblock{ J. Appl. Probability, v. 49, 303-318.}


\end{thebibliography}
\end{document}